%% file: ram_split_model.tex
\documentclass{article}

\usepackage[utf8]{inputenc}
\usepackage{amsmath,amsfonts,amssymb,amsthm,pb-diagram,picinpar,graphicx,color}
\usepackage{algorithm}
\usepackage{amscd}
\usepackage{algorithmic}
\usepackage{colortbl}
\usepackage{enumerate}
\usepackage{hyperref}
\usepackage{bm}
\usepackage{tikz}
\usepackage{pgf}
\usepackage[T1]{fontenc}
\usepackage{enumitem}
\usepackage{comment}
\usetikzlibrary{shapes,matrix,calc,arrows,decorations.markings,knots,patterns,intersections,cd}

\newcommand{\CC}{\mathbb{C}}
\newcommand{\PP}{\mathbb{P}}

\newcommand{\mcO}{\mathcal{O}}
\newcommand{\mcQ}{\mathcal{Q}}
\newcommand{\mcT}{\mathcal{T}}

\newcommand{\ff}{\mathfrak {f}}

\newcommand{\sing}{\mathop{\mathrm{Sing}}\nolimits}
\newcommand{\Red}{\mathop{\mathrm{Red}}\nolimits}

\newcommand{\MW}{\mathop{\mathrm{MW}}\nolimits}

\newcommand{\I}{\mathop {\mathrm {I}}\nolimits}
\newcommand{\II}{\mathop {\mathrm {II}}\nolimits}
\newcommand{\III}{\mathop {\mathrm {III}}\nolimits}
\newcommand{\IV}{\mathop {\mathrm {IV}}\nolimits}

\newtheorem{thm}{Theorem}[section]
\newtheorem{thm0}{Theorem}
\newtheorem{lem}[thm]{Lemma}     

\theoremstyle{definition}
\newtheorem{defin}[thm]{Definition}

\newtheorem{ex}[thm]{Example}

\theoremstyle{remark}
\newtheorem{rem}[thm]{Remark}

\renewcommand{\thesubparagraph}{\theparagraph.\@arabic\c@subparagraph}
\begin{document}
\begin{center}

{\LARGE \bf  Ramified and Split models of elliptic surfaces and bitangent lines of quartic curves}

\bigskip

 {\large \bf Shinzo Bannai\footnote{Partially supported by JSPS KAKENHI Grant Numbers JP18K03263, JP23K03042},  Hiro-o Tokunaga\footnote{Partially supported by JSPS KAKENHI Grant Number JP20K03561} and Emiko Yorisaki}
 
 \end{center}

\begin{abstract}
We define ramified and split models of elliptic surfaces and study  the relation between the two models. We focus on certain rational 
elliptic surfaces from these points of views  and as an application, we give an observation
 on bitantgent lines of an irreducible quartic with at most nodes.\\
{\bf Keywords: } elliptic surfaces, split model, ramified model, bitangent lines of a plane
quartic.\\
{\bf MSC2020: } 14H50, 14H52, 14H27
\end{abstract}

\bigskip

{\LARGE \bf Introduction}

\medskip 

Let $\varphi: S \to C$ be an elliptic surface over a smooth projective
curve $C$ satisfying:
\begin{itemize}
    \item $\varphi$ is relatively minimal,
    \item $\varphi$ has a distinguished section  denoted by $O$, and
    \item $\varphi$ has at least one singular fiber.
\end{itemize}
Throughout this article, we assume that these three conditions are satisfied for elliptic surfaces.
By the second condition, the generic fiber, $E_S$, of $\varphi$ is a curve of genus $1$ over the 
rational function field, $\CC(C)$, of $C$ with a rational point,
which we also denote by $O$. Namely, $(E_S, O)$ is an elliptic curve over $\CC(C)$. We can define an abelian
group structure,  whose  addition we denote by $\dot{+}$, on $E_S$ with $O$ as the zero element in the well-known manner. We denote the inversion
with respect to $\dot{+}$ by $[-1]$. We denote the set of $\CC(C)$-rational points of $E_S$ by $E_S(\CC(C))$. Under our settings, $E_S(\CC(C))$ can be identified with the set of 
sections of $S$, which we denote by $\MW(S)$. For $P \in E_S(\CC(C))$, we denote its
corresponding section by $s_P$. Conversely, for $s \in \MW(S)$, we denote its corresponding
$\CC(C)$-rational point, $P_s$,  given by $s|_{E_S}$.
Given the distinguished section $O$, it is known that $S$ can be canonically realised as the minimal resolution of a double
cover $S'_O$ of a ruled surface $\Sigma$ over $C$,  namely as the  Weierstrass model $\tilde{f}_{O}: S \to \Sigma$ (See \cite[II]{miranda-BTES} and see \cite{horikawa} for the canonical resolution). 
The covering transformation of this double cover is induced by $[-1]$ on the generic fiber.
If $\varphi$ has a section $s$ distinct from $O$, 
we have an involution on $E_S$ given by
$Q \mapsto [-1]Q \dot{+}P_s$, $Q \in E_S$. This induces a fiber
preserving involution on $S$, i.e., $t_{P_s}\circ [-1]$,
$t_{P_s}$ being the translation by $P_s$. We denote  this involution by $\sigma_s$. Conversely any fiber preserving involution which is distinct from
$[-1]$ and has fixed points on a general fiber is obtained this way.
As we see in Section \ref{sec:fixed_points}, $\sigma_s$ has no isolated points.
Hence, by considering the quotient $S/\langle\sigma_s\rangle$ by $\sigma_s$, $S$ is realised as the canonical resolution of a double cover, $S'_{O, P_s}$, of 
a ruled surface $\Sigma'$ over $C$. We denote  the induced morphism by $\tilde {f}_{O, P_s}$. The following is a diagram
for other morphisms related to these two realizations ${\tilde f}_{O}$ and ${\tilde f}_{O, P_s}$:

\begin{center}
\begin{tikzpicture}
\node (1) at (0,0) {$\Sigma$};
\node (2) at (2,0) {$S/[-1]$};
\node (3) at  (4,0) {$S/\langle \sigma_s \rangle$};
\node (4) at (6, 0) {$\Sigma'$};
\node (5) at (0, 2) {$S'_O$};
\node (6) at (3, 2) {$S$};
\node (7) at (6, 2) {$S'_{O,P_s}$};
\draw [->] (6) to node[left]{$f_O$} (2);
\draw [->] (5) to node[left]{$f'_O$} (1);
\draw [->] (6) to node[above]{$\mu_O$} (5);
\draw [->] (6) to node[right]{$f_{O,P_s}$} (3);
\draw [->] (6) to node[above]{$\mu_{O, P_s}$} (7);
\draw [->] (7) to node[right]{$f'_{O,P_s}$} (4);
\draw[->] (3) to node[below]{$q_{O,P_s}$} (4);
\draw  [->] (2) to node[below]{$q_O$} (1);
\draw [->] (6) to node[above]{$\tilde{f}_O$} (1);
\draw [->] (6) to node [above right]{$\tilde{f}_{O,P_s}$} (4);
\end{tikzpicture}
\end{center}
Here $f'_{O} : S_O' \to \Sigma$ (resp. $f'_{O, P_s}: S'_{O,P_s} \to \Sigma'$) is the Stein factorization of ${\tilde f}_{O}$ (resp.
${\tilde f}_{O, P_s}$),  $f_O$ and $f_{O, P_s}$ are the quotient morphisms, 
$\mu_O$ and $\mu_{O, P_s}$ are the canonical resolution and $q_O$ and $q_{O, P_s}$ are  compositions of blowing-ups so that the branch loci of $f_{O}$ and $f_{O, P_s}$
become smooth. Note that ${\tilde f}_O = q_O\circ f_O$
and ${\tilde f}_{O, P_s} = q_{O, P_s}\circ f_{O, P_s}$.

\begin{defin}\label{def:model}{\rm Let $\varphi: S \to C$ be an elliptic surface as above.
\begin{itemize}
 \item We call the realization ${\tilde f}_O : S \to \Sigma$ the ramified model for $S$ with respect to $O$.
    \item We call the realization  ${\tilde f}_{O, P_s} : S \to \Sigma'$ the split model for $S$ with respect to $O, P_s$.   
\end{itemize}
}
\end{defin}
We study the action of $\sigma_s$ on singular
fibers, analyse the branch locus of $\tilde{f}_{O, P_s}: S \to \Sigma'$ and consider how the above two models
are related. Namely, when the generic fiber of $S$ as a ramified model, $(E_S, O)$, is given by an  Weierstrass equation of the form
\[
 y^2 = x^3 + ax^2 + bx + c, \quad a, b, c \in \CC(C),
\]
and $P= (x_P, y_P)$ is also given, we find 
a defining equation of the generic fiber $(E_S, \{O, P_s\})$ of 
the corresponding split model $S'_{O,P_s}$ of the form
\[
{y'}^2 = (x'^2 + a')^2 + b'x' + c'=  {x'}^4 + 2a'{x'}^2 +b'{x'} + a'^2 + c', \quad a', b', c' \in \CC(C).
\]
%
We then consider how
the branch loci of ${\tilde f}_O$ and ${\tilde f}_{O, P_s}$ are related.

In particular, we are interested in the case when $S$ is a
rational elliptic surface and $s\cdot O = 0$ (i.e., $s$ is an integral section, see Section \ref{sec:notation}).
In this case, we have the following facts:
\begin{itemize}
    \item $\Sigma$ is the Hirzebruch surface of degree $2$, $\Sigma_2$, and the branch locus of 
    $f'_{O}$ is of the form $\Delta_0 + \mcT$, where $\Delta_0$ is the section with $\Delta^2_0 = -2$ and $\mcT$ is
    a tri-section such that $\mcT\cap \Delta = \emptyset$. Also $\Delta_0$
    gives rise to the  section $O$.
    \item $\Sigma'$ is the Hirzebruch surface of degree $1$, $\Sigma_1$, i.e.,
    a surface obtained by a blowing up of $\PP^2$ at a point $z_o$. The branch locus of $\tilde{f}_{O, P_s}$ is the preimage of a reduced
    quartic $\mcQ^\prime$ such that $\mcQ^\prime$ is not a union of concurrent $4$ lines and $z_o \not\in \mcQ^\prime$. The $(-1)$ curve on $\Sigma'$ gives
    rise to the sections $O$ and $P_s$. When $P_s$ is clear from the context, we put $O^+=O$ and $O^-=P_s$ for simplicity. 
    Note that
    we also identify $O^{\pm}$ with the rational points determined by
    them.
\end{itemize}


 We apply the above settings  to the study of  quratics and their bitangent lines in the following way: Let $\mcQ$ be a quartic not being concurrent $4$ lines. Choose a point $z_o \in \PP^2$. We can attach a rational elliptic surface $S_{\mcQ, z_o}$ with $(\mcQ, z_o)$. As we see in  \cite{shioda-E7}, \cite{bannai-tokunaga15, bannai-tokunaga17},\cite{tokunaga10}, 
 if $z_o \in \mcQ$, $S_{\mcQ, z_o}$ can be considered as 
 as a ramified model.  If $z_o \not\in \mcQ$,   
 $S_{\mcQ, z_o}$ is regarded as a split model as in  \cite{kuwata2005}. 
 Our viewpoint in this article is that  we consider both ramified and 
 split realizations for a given rational elliptic surface at one time.
 More prescisely, we consider a ramified model $S_{\mcQ, z_o}$ and a split model $S_{\mcQ', z_o'}$ for a given rational 
 elliptic surface $\varphi: S \to \PP^1$ at one time.
 For both cases, $\varphi$ is  induced by the
 pencil of lines $\Lambda_{z_o}$ (resp. $\Lambda_{z'_o}$) through $z_o$ (resp. $z'_o$)
 and singular fibers  arise from members of 
 $\Lambda_{z_o}$ (resp. $\Lambda_{z'_o}$) meeting $\mcQ$ (resp. $\mcQ'$) in a special way, i.e., tangent to $\mcQ$ (resp. $\mcQ'$), passing through  a singularity of $\mcQ$
(resp. $\mcQ'$)  and  so on.

  Choose $z_o \in \mcQ$ and let $\varphi_{\mcQ,z_o} : S_{\mcQ, z_o} \to \PP^1$ be the ramified model attached
  to $(\mcQ, z_o)$. Choose a section $s$ of $\varphi_{\mcQ, z_o}$ and consider the split model 
  $\tilde{f}_{O, P_s}: S_{\mcQ, z_o} \to \Sigma_1$. Then we obtain a quartic $\mcQ'$ and a point $z'_o \in \PP^2\setminus \mcQ'$
  such that {the split model $S_{\mcQ', z'_o}$ coincides with $S_{\mcQ,z_o}$.}
  For each singular fiber $F$, we consider
  relations between the type of $F$ and at which component $s$ meets with $F$. This makes us possible to determine
  how the line that is the image of $F$ meets $\mcQ'$.
We apply this consideration to count the number of concurrent bitangent lines for  an irreducible quartic $\mcQ$
with at most nodes. Here a line $l$ is said to be a bitangent to $\mcQ$ if either $l$ is tangent to $\mcQ$ at two distinct smooth points or $l$ is a $4$-fold tangent  at a smooth point of $\mcQ$.  The former case is said to be ordinary, while the latter is 
said to be special. Now our statement is as follows:

\begin{thm0}\label{thm:concurrent} Let $\mcQ$ be an irreducible quartic with at most nodes. Then we have

\begin{center}
\begin{tabular}{|c|c|}\hline
The number of nodes & the number of concurrent bitangent lines \\ \hline
 $0$  & $\le 4$ \\
 $1$  & $\le 4$ \\
 $2$  & $\le 4$ \\
 $3$  & $\le 3$ \\ \hline
\end{tabular}
\end{center}
All maximal cases occur. Moreover, for each maximal case, if $\mcQ$ has $k$ ordinary bitangents and $l$ 
special bitangents, then the  only following cases occur:
\begin{center}
\begin{tabular}{|c|c|}\hline
The number of nodes & Possible pairs $(k,l)$ \\ \hline
 $0$  & $(4,0), (3,1), (2,2), (0, 4)$ \\
 $1$  & $(4,0), (3, 1), (2, 2)$ \\
 $2$  & $(4, 0)$ \\
 $3$  & $(3,0)$ \\ \hline
\end{tabular}
\end{center}
\end{thm0}

This paper is organized as follows: In Sections \ref{sec:notation} and \ref{sec:transformation},
we explain terminologies and transformations between ramified and 
split models, respectively. We study fixed points of $\sigma_s$ in
Section \ref{sec:fixed_points} and actions of $\sigma_s$ on sigular fibers in Section \ref{sec:action}.
We prove Theorem~\ref{thm:concurrent} in Section \ref{sec:proof}. Explicit
examples are given in Section \ref{sec:examples}.
\input preliminaries.tex


\input transformation.tex

\input fixed_locus.tex

\input singular-fibers.tex

\input proof_of_theorem1.tex

\input example.tex

\bibliographystyle{spmpsci}
\bibliography{biblio.bib}

\noindent Shinzo BANNAI \\
Department of Applied Mathematics, Faculty of Science, \\
Okayama University of Science, 1-1 Ridai-cho, Kita-ku,
Okayama 700-0005 JAPAN \\
{\tt bannai@ous.ac.jp}\\

\noindent  Hiro-o TOKUNAGA and Emiko YORISAKI\\
Department of Mathematical  Sciences, Graduate School of Science, \\
Tokyo Metropolitan University, 1-1 Minami-Ohsawa, Hachiohji 192-0397 JAPAN \\
{\tt tokunaga@tmu.ac.jp}
\end{document}

%% file: preliminaries.tex

\section{Notation and terminologies}\label{sec:notation}

As for theory of elliptic surfaces and their related results which we need
in this article, our main references are
\cite{bannai-tokunaga15,
kodaira, kuwata2005, miranda-BTES, shioda90, tokunaga12}. 

We keep the notation and terminology as in the Introduction. As for singular fibers, we use Kodaira's
notation (\cite{kodaira}) to describe their types and label irreducible components as in \cite{tokunaga12,tokunaga14}. For $\varphi : S \to C$, we set $\sing(\varphi):= \{v \in C \mid \mbox{$F_v = \varphi^{-1}(v)$ is singular.}\}$ and $\Red(\varphi):= \{v\in \sing(\varphi)\mid \mbox{$F_v$ is reducible.}\}$. {We write the irreducible decomposition of a singular fiber $F_v$ by 
$F_v = \Theta_{v, 0} + \sum_{i}a_{v,i}\Theta_{v,i}$.
 Here  $\Theta_{v,0}$ is
the unique component with $O\cdot \Theta_{v,0} = 1$ and  
we  call $\Theta_{v,0}$ the identity component.}
For a singular fiber $F_v$, put $F_v^{\sharp} = \bigcup_{a_{v,i} = 1}\Theta_{v, i}^{\sharp},\, 
\Theta_{v,i}^{\sharp}:= \Theta_{v,i}\setminus(\mbox{singular points of $(F_v)_{\mathrm{red}}$})$. Here we put 
$a_{v, 0} = 1$. By \cite[\S 9]{kodaira}, $F_v^{\sharp}$ has a structure of ablelian group. We define
$G_{F_v^{\sharp}}$ to be a finite group as in \cite[p.82]{tokunaga12}. Put $G_{\sing(\varphi)}:= \oplus_{v\in \sing(\varphi)} G_{F_v^{\sharp}}$. As in \cite[p.83]{tokunaga12}, we define the homomorphism $\gamma: \MW(S) \to G_{\sing(\varphi)}$. By abuse of noation, we write $\gamma(P) = \gamma(s_P)$ for $P\in E_S(\CC(C))$. 
$\gamma(s)$ describes at which irreducible component $s$ meets at $F_v$.
We denote the height pairing define by Shioda in \cite{shioda90} by $\langle \bullet, \bullet\rangle$.
We use the notation, terminologies and properties in \cite{shioda90} freely for $\langle\bullet, \bullet\rangle$.



%% file: transformation.tex
\section{Transformations between  a ramified model and a split 
model}\label{sec:transformation}

Let $K$ be a perfect field of ${\rm char}(K)\neq 2$. Let $(E, \{O^+, O^-\})$ be an elliptic curve 
given by 
\[
y'^2 =  (x'^2 + a')^2 + b'x' + c', \quad  a', b', c' \in K.
\]

A transformation from this
equation to that of a ramified model $(E, O^+)$ is well-known, for example, it is given as follows (cf. \cite[p. 35]{cassels}):
\begin{equation}
\label{eq:transformation1}
x = y' + x'^2 + a', \quad y = \sqrt{2}x' x + \frac {b'}{2\sqrt{2}}.
\end{equation}

In \cite[p.77, p139]{mordell} and \cite[\S 2]{shioda-Weierstrass}, similar arguments are found. With the transformation (\ref{eq:transformation1}), we have an equation of
$(E, O^+)$:
\begin{equation}
\label{eq:ramified0}
y^2 = x^3 - 2a'x^2 - c'x + \frac{b'^2}8
\end{equation}
\begin{rem}{\rm The above $x$ and $y$ give rational functions such that
poles of $x$ and $y$ are $2O^+$ and $3O^+$, respectively.
}
\end{rem}

In this section, we give an explicit description for the converse of the above transformation.
More precisely, given a section $s$ of $S$, we consider a transformation from $(E_S, O)$
to $(E_S, \{O, P_s\})$ in terms of equations of generic fibers. 

Let $\varphi:S \rightarrow C$ be an elliptic surface as in the Introduction.
Let $(E_S, O)$ be the generic fiber of a ramified model for $S$ given by
\begin{equation}
\label{eq:ramified}
y^2 = x^3 + ax^2 + bx + c
\end{equation}
where $a,b,c \in \mathbb{C}(C)$.
We choose a $\mathbb{C}(C)$-rational point on $(E_S, O)$ as $P:=(x_P,y_P)$.

Consider a transformation $w_P$ given by
\begin{equation}
\label{eq:transformation2}
x = y'+x'^2-\frac{1}{2}(x_P + a) \qquad
y = \sqrt{2}x'(x-x_P) -y_P.
\end{equation}
Since $x_P, y_P$ satisfy the equation~\ref{eq:ramified}, we 
have an equation of the split model $(E_S, (O, P))$:
\begin{eqnarray*}
(E_S, \{O, P\}):y'^2 
&=& \left( x'^2 -\frac{1}{2}(3x_P + a) \right)^2 - 2\sqrt{2}y_Px' - (3x_P^2 +2ax_P + b).
\end{eqnarray*}
The involution $(x', y') \mapsto (x', -y')$ on $(E_S, \{O, P\})$ is described on $(E_S, O)$ as follows:
As in the Introduction, we regard $[-1]$ as a involution on ramified model $(E_S, O)$, i.e., $[-1]Q = (x_Q,-y_Q)$ for $Q:=(x_Q,y_Q) \in (E_S, O)$.
Let $Q:=(x_Q', y_Q')$ be the corresponding  point to $Q$ on $(E_S, \{O, P\})$.
Then we have
\begin{eqnarray*}
 \left(x_Q',y_Q' \right) &\mapsto& \left(x_Q,y_Q \right) \\
&=& \left( y_Q'+x_Q'^2-\frac{1}{2}(x_P+a),\sqrt{2}x_Q'(x_Q-x_P)- y_P \right)  \\
 \left(x_Q',-y_Q' \right) &\mapsto& \left( -y_Q'+x_Q'^2-\frac{1}{2}(x_P+a),\sqrt{2}x_Q' \left(-y_Q'+x_Q'^2-\frac{1}{2}(x_P+a)-x_P \right) - y_P \right) \\
&=& \left( \left(y_Q'+x_Q'^2-\frac{1}{2}(x_P+a) \right)-2y_Q', \sqrt{2}x_Q' \left(-y_Q'+x_Q'^2-\frac{1}{2}(x_P+a)-x_P \right) -y_P \right) \\
&=& \left(x_Q-2y_Q', \sqrt{2}x_Q'\left( \left(x_Q-2y_Q' \right)-x_P \right) - y_P \right) \\
&=& \left(x_Q-2y_Q', \left(\sqrt{2}x_Q' \left(x_Q-x_P \right) - y_P \right) - 2\sqrt{2}x_Q'y_Q' \right) \\
&=& \left(x_Q-2y_Q', y_Q-2\sqrt{2}x_Q'y_Q' \right)
\end{eqnarray*}

Now by explicit computation for $(x_Q, -y_Q)\dot{+}(x_P, y_P)$ we see that
\[
 \left(x_Q-2y_Q', y_Q-2\sqrt{2}x_Q'y_Q' \right) =
(x_Q,-y_Q) \dot{+} (x_P, y_P) = [-1]Q \dot{+} P.
\]
From the above observation, the map $(x_Q',y_Q') \mapsto (x_Q', -y_Q')$ defines a involution $(E_S, \{O, P\})$ which coincide the involution $t_{P}\circ [-1]$ on $(E_S, O)$.
If $Q$ is a fixed point of $t_P\circ [-1]$ and $P$ is 2-torsion on $(E_S, O)$, then $Q \dot{+} P$ also becomes a fixed point of $t_P\circ[-1]$.
Indeed, suppose $Q \in E_S$ is a fixed point and $2P=O$, we have

\[
t_P\circ[-1](Q \dot{+} P) = [-1](Q \dot{+} P) \dot{+} P = [-1]Q = [-1]([-1]Q \dot{+} P)=  Q \dot{+}[-1]P =  Q \dot{+}P.
\]

%% file: fixed_locus.tex

\section{The set of fixed points of a fiber preserving involution \texorpdfstring{$\sigma$}{}}
\label{sec:fixed_points}

Let $\varphi: S \to C$ be an elliptic surface and let $\sigma$ be a fiber preserving involution on $S$ which has fixed points on a general fiber. 
We denote the set of fixed points of $\sigma$
by $S^{\sigma}$. By considering the action of $\sigma$ on the tangent space at each fixed point, we see that
$S^{\sigma}$ is a union of isolated fixed points and some disjoint smooth curves on $S$. Write
\[
S^{\sigma} = \{p_1, \ldots, p_{n_1}\}\cup(\cup_{k=1}^{n_2}D_k),
\]
where $p_j$'s are isolated points and $D_k$'s denote a smooth curve.

\begin{lem}\label{lem:fixed-points}{$S^{\sigma}$ has no isolated points, i.e., $n_1 = 0$}
\end{lem}
\proof By the holomorphic Lefschetz Theorem and the arguments in  \cite[p.567-568]{atiyah-singer}, we have
\[
(\ast) \quad \quad  \sum_{p = 0}^2 (-1)^p {\mathrm{trace}}(\sigma^*|_{H^p(S, \mcO_S)} ) = \frac {n_1}4 + \sum_{k=1}^{n_2} 
\left (\frac {1 - g(D_k)}2 + \frac 14D_k^2\right ),
\]
where $g(D_k)$ is the genus of $D_k$. Since 
$H^1(S, \mcO_S) = \overline H^0(S, \Omega_S^1) = \overline{\varphi^*H^0(C, \Omega_C^1)}$, $\sigma^*$ is trivial on $H^1(S, \mcO_S)$. Hence ${\mathrm{trace}}(\sigma^*|_{H^1(S, \mcO_S)} ) = q(S)$. Since $S/\langle \sigma \rangle$ has
at most quotient singularities and its smooth model is birationally 
equivalent to a ruled surface,
$H^0(S/\langle \sigma \rangle, \Omega^2_S) = (H^0(S, \Omega^2_S))^{\sigma} = \{0\}$
by \cite[Proposition 9.24]{ueno}. Hence 
${\mathrm{trace}}(\sigma^*|_{H^2(S, \mcO_S)} ) = - p_g$. Therefore the left hand side of $(\ast)$ is $1 - q - p_g$.
As for the right hand side, by the adjunction formula we have
\begin{eqnarray*}
 &  & \frac {1 - g(D_k)}2 + \frac 14 D^2_k \\
 & = & \frac 14( -D_k^2 -D_k\cdot K_S) + D^2_k) = - \frac14D_k\cdot K_S.
 \end{eqnarray*}
 
 Thus we have
 \[
 \sum_{k=1}^{n_2} 
\left (\frac {1 - g(D_k)}2 + \frac 14D_k^2\right ) = -\frac14\left (\sum_{k=1}^{n_2}D_k\right )\cdot K_S.
\]
By the canonical bundle formula (see \cite[Ch. V, \S 12]{bhpv},  
\[
K_S \approx (\chi(\mcO_S) + 2g(C) - 2)F = (\chi(\mcO_S) + 2q - 2)F.
\]
This implies the right hand side of $(\ast)$ is
\[
 -\frac14\left (\sum_{k=1}^{n_2}D_k\right )\cdot K_S = - \frac 14 (\chi(\mcO_S) + 2q - 2)\sum_{k=1}^{n_2}D_k\cdot F.
 \]
 As $\sigma$ has a fixed point on a general fiber $F$, 
 $F/\langle\sigma\rangle = \PP^1$. Hence $\sum_{k=1}^{n_2}D_kF =4$. This implies that $n_1 = 0$, i.e., $S^{\sigma}$ has
 no isolated fixed point. \qed

\begin{rem}{In \cite[Lemma 4.1]{oguiso-zhang}, the holomorphic  Lefschetz Theorem was used in order to 
consider the fixed loci of automorphisms of oder five on K3 surfaces. Our argument above can be a simplified 
version.}
\end{rem}

%% file: singular-fibers.tex

\section{Actions of \texorpdfstring{$\sigma$}{Lg} on singular fibers}\label{sec:action}

Let $\varphi : S \to C$ and $\sigma$ be as in Section~\ref{sec:fixed_points}. For simplicity, 
In this section, we consider how $\sigma$ acts on singular fibers of $S$. Let us start with the following lemma.

\begin{lem}\label{lem:fixed-points2}{Let $D_1$ and $D_2$ be irreducible curves on $S$
 such that {\rm (i)} $\sigma(D_1\cup D_2) =
D_1\cup D_2$, {\rm (ii)} $D_1 \cap D_2 = \{p\}$ and both of $D_i, \, (i = 1, 2)$ are smooth at $p$. Then one of the following
 holds:
 
  \begin{enumerate}
  \item[\rm (a)] $D_i \not\subset S^{\sigma}, i =1, 2$ and there exists an irreducible component $D$ of $S^{\sigma}$ such that $p \in D$.
 
 \item[\rm (b)] Either $D_1\subset S^{\sigma}, D_2\not\subset S^{\sigma}$ or 
 $D_2\subset S^{\sigma}, D_1\not\subset S^{\sigma}$ holds.
\end{enumerate}
}
\end{lem}

\proof By our condition,  $p \in S^{\sigma}$. If $D_1, D_2 \not\subset S^{\sigma}$,there exists an irreducible component through $p$  by
 Lemma~\ref{lem:fixed-points}. 
 Assume that $S^{\sigma}$ contains at least one of $D_1$ and $D_2$. 
 If $D_1$ and $D_2$ meet at $p$ transversely, by considering the action of $\sigma$ on the
 tangent space at $p$, we see that either (b) follows. 
 Now we may assume that $D_1 \subset S^{\sigma}$.
 If $D_1$ and $D_2$ does not meet at $p$ transversely, 
 we can choose local coordinates $(z_1, z_2)$ such that $D_1$ and $D_2$ are locally given by $z_2 = 0$ and $z_2 - \phi(z_1) = 0$ respectively 
 around $p =(0,0)$
 and
 $\sigma$ is given by
 $(z_1, z_2) \mapsto (z_1, -z_2)$, since both of $D_i$ $(i = 1, 2)$ are smooth. Hence $\sigma(D_2) \neq
 D_1, D_2$, which contradicts to the assumption (i). \qed

 \medskip 
 Before we go on to look into actions of singular fibers in detail, we recall the following facts:
 $\sigma$ induces an involution of the generic fiber $E_S$, which we also denote by $\sigma$ for simplicity, and  $\sigma$
 is a composition of the inversion and a translation. 
 Put $s_{\sigma} = \sigma(O)$. Then we see that
 $\sigma$ is given by 
 $\sigma(Q) = [-1]Q\dot{+}P_{s_{\sigma}}, Q 
\in E_S(\CC(C))$. 
 We also label the irreducible components of singular fibers for each type as in \cite{tokunaga12, tokunaga14}. Under these
 settings, by \cite[\S 9]{kodaira},   $\sigma$ acts irreducible components of reducible singular fibers as follows:
 
 Here the middle column shows at which component $s_{\sigma}$ meets and the right column shows how each
 irreducible components are mapped;
 \begin{center}
 \begin{tabular}{|c|c|c|}\hline 
Type of singular fiber & $s_\sigma\cdot\Theta_\bullet=1$ & action of $\sigma$\\
 \hline
 Type $\I_b$ $(b \ge 2)$ & $\Theta_l$ & $\Theta_k \mapsto \Theta_{b - k + l \bmod b}.$\\ \hline
 
 Type $\I^*_b$ ($b \ge 0$, even) &$\Theta_0$ & 
 $\displaystyle{\begin{array}{c} 
 \Theta_k \mapsto \Theta_k,  \, k = 0, 4, \ldots, b+ 4,  \\
 \Theta_{ij} \mapsto \Theta_{ij},  \, (i, j) = (1, 0), (0, 1), (1, 1). \end{array}}$ \\ \hline
 
 Type $\I^*_b$ ($b \ge 0$, even) &$\Theta_{10}$ & 
 $\displaystyle{\begin{array}{c} 
 \Theta_k \mapsto \Theta_k,  \, k =  4, \ldots, b+ 4,  \\
 \Theta_{0} \mapsto \Theta_{10} \mapsto \Theta_0, \,
 \Theta_{01} \mapsto \Theta_{11} \mapsto \Theta_{01}. \end{array}}$ \\ \hline
 
  Type $\I^*_b$ ($b \ge 0$, even) &$\Theta_{01}$ & 
 $\displaystyle{\begin{array}{c} 
 \Theta_k \mapsto \Theta_{b+8 -k},  \, k =  4, \ldots, b+ 4,  \\
 \Theta_{0} \mapsto \Theta_{01} \mapsto \Theta_0, \, 
 \Theta_{10} \mapsto \Theta_{11} \mapsto \Theta_{10}. \end{array}}$ \\ \hline
 
  Type $\I^*_b$ ($b \ge 0$, even) &$\Theta_{11}$ & 
 $\displaystyle{\begin{array}{c} 
 \Theta_k \mapsto \Theta_{b + 8-k},  \, k =  4, \ldots, b+ 4, \\
 \Theta_{0} \mapsto \Theta_{11} \mapsto \Theta_0, \, 
 \Theta_{10} \mapsto \Theta_{01} \mapsto \Theta_{10}. \end{array}}$ \\ \hline
 
Type $\I^*_b$ ($b \ge 	1$, odd)  &$\Theta_{0}$ &  $\Theta_k \mapsto \Theta_k \, (\forall k)$. \\ \hline

 Type $\I^*_b$ ($b \ge 1$, odd) &$\Theta_{1}$ & 
 $\displaystyle{\begin{array}{c} 
 \Theta_k \mapsto \Theta_{b + 8-k},  \, k =  4, \ldots, b+ 4, \\
 \Theta_{0} \mapsto \Theta_{1} \mapsto \Theta_0, \,
 \Theta_{2} \mapsto \Theta_{3} \mapsto \Theta_{2}. \end{array}}$ \\ \hline
 
 Type $\I^*_b$ ($b \ge 1$, odd) &$\Theta_{2}$ & 
 $\displaystyle{\begin{array}{c} 
 \Theta_k \mapsto \Theta_{k},  \, k =  4, \ldots, b+ 4, \\
 \Theta_{0} \mapsto \Theta_{2} \mapsto \Theta_0, \,
 \Theta_{1} \mapsto \Theta_{3} \mapsto \Theta_{1}. \end{array}}$ \\ \hline

 Type $\I^*_b$ ($b \ge 1$, odd) &$\Theta_{3}$ & 
 $\displaystyle{\begin{array}{c} 
 \Theta_k \mapsto \Theta_{b+ 8-k},  \, k =  4, \ldots, b+ 4, \\
 \Theta_{0} \mapsto \Theta_{3} \mapsto \Theta_0, \, 
 \Theta_{1} \mapsto \Theta_{2} \mapsto \Theta_{1}. \end{array}}$\\ \hline
 
Type $\II^*$ & $\Theta_0$ & $\Theta_k \mapsto \Theta_k \, (\forall k)$ \\ \hline

Type $\III$ & $\Theta_0$ & $\Theta_k \mapsto \Theta_k \, k = 0,1$. \\ \hline
Type $\III$ & $\Theta_1$ & $\Theta_0 \mapsto \Theta_1 \mapsto \Theta_0$. \\ \hline
Type $\III^*$ & $\Theta_0$ & $\Theta_k \mapsto \Theta_k, \, (\forall k)$. \\ \hline
Type $\III^*$ &$\Theta_{1}$ & 
 $\displaystyle{\begin{array}{c} 
 \Theta_{0} \mapsto \Theta_{1} \mapsto \Theta_0,\\
 \Theta_k \mapsto \Theta_{9-k} \mapsto \Theta_{k},  \, k = 2, 3, \\
  \Theta_k \mapsto \Theta_{k},  \, k = 4, 5  \end{array}}$ \\ \hline

Type $\IV$ & $\Theta_0$ & $\Theta_0 \mapsto \Theta_0, \Theta_1 \mapsto \Theta_2 \mapsto \Theta_1$ \\ \hline

Type $\IV$ & $\Theta_1$ & $\Theta_0 \mapsto \Theta_1 \mapsto \Theta_0, \Theta_2 \mapsto \Theta_2$ \\ \hline

Type $\IV$ & $\Theta_2$ & $\Theta_0 \mapsto \Theta_2 \mapsto \Theta_0, \Theta_1 \mapsto \Theta_1$ \\ \hline
Type $\IV^*$ &$\Theta_{0}$ & 
 $\displaystyle{\begin{array}{c} 
 \Theta_k \mapsto \Theta_{k}  \, k = 0, 3, 6 \\
  \Theta_1 \mapsto \Theta_{2} \mapsto \Theta_1, \, 
   \Theta_3 \mapsto \Theta_4 \mapsto \Theta_3  \end{array}}$\\ \hline
 Type $\IV^*$ &$\Theta_{1}$ & 
 $\displaystyle{\begin{array}{c} 
  \Theta_k \mapsto \Theta_{k}  \, k = 2, 5, 6 \\
  \Theta_0 \mapsto \Theta_{1} \mapsto \Theta_0, \, 
   \Theta_3 \mapsto \Theta_4 \mapsto \Theta_3  \end{array}}$ \\ \hline
  Type $\IV^*$ &$\Theta_{2}$ & 
 $\displaystyle{\begin{array}{c} 
  \Theta_k \mapsto \Theta_{k}  \, k = 1, 5, 6 \\
  \Theta_0 \mapsto \Theta_{2} \mapsto \Theta_0, \,
  \Theta_3 \mapsto \Theta_5 \mapsto \Theta_3  \end{array}}$ \\ \hline
 
 \end{tabular}
 \end{center}
 
 We next describe the singularities of the branch locus $\mcQ^\prime$ of $f^\prime_{O, P_s}$ and its relation with the fibers of $\Sigma^\prime$.  Let $F$ denote a singular fiber of $S$ and let $F^\prime$ be the image of $F$ in $\Sigma^\prime$. Furthermore, let $\overline{\Theta}_{\bullet}=f_{O, P_{s_{\sigma}}}
(\Theta_{\bullet})$ in $S/\langle\sigma\rangle$. Note that since we are assuming that $O^+=O$ and $O^-=P_{s_{\sigma}}$,  $\Theta_0$ maps surjectively onto $F^\prime$. 

\begin{itemize}
    \item $F$ is of Type $\I_{2n}$ and $l$ is even.
    
    Suppose, $s_\sigma\cdot\Theta_l=1$ and $l=2l^\prime$ for some $0\leqq l^\prime<n$. Then, since $\sigma(\Theta_k)=\Theta_{2n-k+l \mod 2n}$ the component $\Theta_{k}$ is fixed by $\sigma$ if $k-l^\prime \equiv 0 \mod{n}$. Hence, there are two components $\Theta_{l^\prime}$ and $\Theta_{l^\prime+n}$ fixed by $\sigma$. The chain $\Theta_{l^\prime+1}, \ldots, \Theta_{l^\prime+(n-1)}$ is transformed to the chain $\Theta_{l^\prime-1 \mod 2n}, \ldots, \Theta_{l^\prime-(n-1)\mod 2n}$, respectively so $\Theta_{l^\prime}$ and $\Theta_{l^\prime+n}$ are not point-wise fixed, and have two fixed points each, which are the intersection points with the ramification locus.
     Then
    $\overline{\Theta}_{l^\prime+i}=\overline{\Theta}_{l^\prime-i \mod 2n}$ ($i=1,\ldots, n-1$) and  $\overline{\Theta}_{l^\prime}, \overline{\Theta}_{l^\prime+1}\ldots, \overline{\Theta}_{l^\prime+n}$ form a single chain with 
    \begin{gather*}
        \overline{\Theta}_{l^\prime+i}\cdot\overline{\Theta}_{l^\prime+i+1}=1, \quad (i=0, \ldots, n-1)\\
        \overline{\Theta}_{l^\prime}^2=\overline{\Theta}_{l^\prime+n}^2=-1\\
        \overline{\Theta}_{l^\prime+i}^2=-2,\quad (i=1, \ldots, n-1).
    \end{gather*}
    Furthermore,  $\overline{\Theta}_{l^\prime},  \overline{\Theta}_{l^\prime+n}$ each intersect with the branch locus transversally at two points.
    Finally, by blowing down $\overline{\Theta}_{l^\prime}, \overline{\Theta}_{l^\prime-1}\ldots, \overline{\Theta}_{1}$ and $\overline{\Theta}_{l^\prime+n}, \overline{\Theta}_{l^\prime+n+1}, \ldots, \overline{\Theta}_{2n-1}$ consecutively in this order, we see that the branch locus has two singular points, one of type $A_{2l^\prime-1}$ and another of type $A_{2(n-l^\prime)-1}$ on  $F^\prime$. Here, the type $A_{-1}$ will denote the situation where the image of the branch locus is smooth and the fiber $F^\prime$ is transversal to the branch locus at two distinct points.
    
    \item $F$ is of Type $\I_{2n}$ and $l$ is odd.
    
    Suppose, $s_\sigma\cdot\Theta_l=1$ and $l=2l^\prime+1$ for some $0\leqq l^\prime<n-1$. In this case, since $k\not\equiv 2n-k+l \mod 2n$ for all $k$, $\sigma$ does not fix any component, of $F$. Furthermore, since $\sigma(\Theta_{l^\prime})=\Theta_{l^\prime+1}$, $\sigma(\Theta_{l^\prime+n})=\Theta_{l^\prime+n+1}$ the points $\{p_1\}=\Theta_{l^\prime}\cap \Theta_{l^\prime+1}$ and $\{p_2\}=\Theta_{l^\prime+n}\cap \Theta_{l^\prime+n+1}$ are fixed by $\sigma$. 
    Hence, $F$ intersects with the ramification locus at these two points. In this case the chain $\Theta_{l^\prime+1}\ldots, \Theta_{l^\prime+n}$ is transformed to the chain $\Theta_{l^\prime}, \Theta_{l^\prime-1\mod 2n}, \ldots, \Theta_{l^\prime-(n-1) \mod 2n}$. Now, 
    $\overline{\Theta}_{l^\prime+i}=\overline{\Theta}_{l^\prime-(i-1)}$ $(i=1, \ldots, n)$ and $\overline{\Theta}_{l^\prime+1}, \overline{\Theta}_{l^\prime+2}\ldots, \overline{\Theta}_{l^\prime+n}$ form a single chain with 
     \begin{gather*}
        \overline{\Theta}_{l^\prime+i}\cdot\overline{\Theta}_{l^\prime+i+1}=1, \quad (i=1, \ldots, n-1)\\
        \overline{\Theta}_{l^\prime+1}^2=\overline{\Theta}_{l^\prime+n}^2=-1\\
        \overline{\Theta}_{l^\prime+i}^2=-2,\quad (i=2, \ldots, n-1).      
    \end{gather*}
    Furthermore, $\overline{\Theta}_{l^\prime}$ and $\overline{\Theta}_{l^\prime+n+1}$ are tangent to the branch locus of $g$. Finally by blowing-down $\overline{\Theta}_{l^\prime}, \overline{\Theta}_{l^\prime-1}\ldots, \overline{\Theta}_{1}$ and  $\overline{\Theta}_{l^\prime+n+1}, \overline{\Theta}_{l^\prime+n+2}\ldots, \overline{\Theta}_{2n-1}$ consecutively in this order, we see that $\mcQ^\prime$ has two singular points,  one of type $A_{2l^\prime}$ and another of type $A_{2(n-l^\prime)-2}$ on  $F^\prime$. Here, the type $A_0$ denotes the situation where $\mcQ^\prime$ and $F^\prime$ are tangent at a smooth point of $\mcQ^\prime$.

    
    \item $F$ is of type $\III$.
    
    The case where $s_\sigma\cdot\Theta_0=1$ is the same as the ramified case, hence $\mcQ^\prime$ has a $A_1$ singularity  and $F^\prime$ is tangent to one of the branches. 
    When $s_\sigma\cdot\Theta_1=1$, $\sigma(\Theta_0)=\Theta_1$ and $\sigma(\Theta_1)=\Theta_0$ and $\{p\}=\Theta_0\cap\Theta_1$ is the unique fixed point on $F$. Since the intersection number of $F$ and the ramification locus is $4$, the ramification locus is tangent to $\Theta_0, \Theta_1$ at $p$. Therefore  $F^\prime=\overline{\Theta}_0=\overline{\Theta}_1$ and $F^\prime$ is tangent to the branch locus  with multiplicity 4.  
\end{itemize}

The cases of the other types of singular fibers can be computed in a similar manner. Note that if $s_\sigma\cdot\Theta_0=1$, then the results will coincide with that of the ramified model. We give a table of the types of singular fibers and the singularities of the branch locus.  The first column describes the type of the singular fiber $F$. The second column describes the component that $s_\sigma$ meets, and the third column describes the singularities of $\mcQ^\prime$ lying on $F^\prime$ and the relation between $\mcQ^\prime$ and $F^\prime$ in $\Sigma^\prime$.
In the case of singular fibers of type $\I_b$, $A_0$ will denote the situation where the branch locus is tangent to  $F$ at a smooth point, and $A_{-1}$ will denote the situation where the branch locus and  $F$ has two distinct transversal intersections.
 
\begin{center}
\begin{tabular}{|c|c|c|}\hline 
Type of singular fiber & $s_\sigma\cdot\Theta_\bullet=1$ & singularities of  $\mcQ^\prime$ and relation with $F^\prime$\\

\hline
 Type $\I_b$ $(b=2n)$ & \begin{tabular}{c}
                        $\Theta_l$, \\
                        $l=2l^\prime$,\\
                        $0\leqq l^\prime \leqq n-1$
                        \end{tabular} & $A_{2l^\prime-1}+A_{2(n-l^\prime)-1}$\\

\hline
  
Type $\I_b$ $(b=2n)$ & \begin{tabular}{c} 
                        $\Theta_l$, \\
                        $l=2l^\prime+1$,  \\
                        $0\leqq l^\prime \leqq n-1$
                        \end{tabular}& $A_{2l^\prime}+A_{2(n-l^\prime)-2}$\\
 
\hline
  
Type $\I_b$ $(b=2n+1)$ & \begin{tabular}{c}
                         $\Theta_{\pm l \mod b}$ \\
                         $0\leqq l \leqq n$
                         \end{tabular} & $A_{2l^\prime-1}+A_{2(n-l^\prime)}$\\
 
%
%
%
 
  \hline
 
 Type $\I^*_b$ ($b=2n$) &$\Theta_0, \Theta_{10}$ & $D_{b+4}$ \\ \hline
 

  Type $\I^*_b$ ($b=2n$) &$\Theta_{01}, \Theta_{11}$ & 
$A_{2n+3}$, $F^\prime$ is contained in tangent cone \\ \hline
 
%
%
Type $\I^*_b$ ($b=2n+1$)  &$\Theta_{0}, \Theta_{2}$ &  $D_{b+4}$ \\ \hline

 Type $\I^*_b$ ($b=2n+1$) &$\Theta_{1}, \Theta_{3}$ & 
 
$A_{2n+4}$, $F$ is contained in tangent cone
 
 \\ \hline

%
%
 
Type $\II^*$ & $\Theta_0$ & $E_8$ \\ \hline

Type $\III$ & $\Theta_0$ & $A_1$, $F^\prime$ is  tangent to a branch with multiplicity 2  \\ \hline

Type $\III$ & $\Theta_1$ & smooth, $F^\prime$ is tangent to $B$ with multiplicity 4\\ \hline
Type $\III^*$ & $\Theta_0$ & $E_7$ \\ \hline
Type $\III^*$ &$\Theta_{1}$ & 
$E_6$, $F^\prime$ is contained in tangent cone \\ \hline

Type $\IV$ & $\Theta_0$ &  $A_2$, $F^\prime$ is contained in tangent cone\\ \hline

Type $\IV$ & $\Theta_1, \Theta_2$ & $A_1$, $F^\prime$ is tangent to a branch with multiplicity 3 \\ \hline


Type $\IV^*$ &$\Theta_{0}$ & 
$E_6$ \\

 \hline

 Type $\IV^*$ &$\Theta_{1}, \Theta_{2}$ & 
$D_4$, $F^\prime$ is contained in tangent cone \\ \hline
%
 \end{tabular}
 \end{center}

%% file: proof_of_theorem1.tex
\section{Proof of Theorem~\ref{thm:concurrent}}\label{sec:proof}

In this section, we restrict ourselves to consider the case of $C = \PP^1$. 

\subsection{The case of  \texorpdfstring{$C=\PP^1$}{}}

When $C = \PP^1$,
$S/[-1]$ can be blow down to the Hirzebruch surface $\Sigma_d$ of degree $d = 2\chi(\mcO_S)$. 
More precisely we realise $S$ as the mimimum resolution of a double cover $f'_O: S'_O \to \Sigma_d$
in the following way:
Let $\Delta_0$ be the section of $\Sigma_d$ with $\Delta_0^2=-d$ and let $\Delta$ be a section linear equivalent
to $\Delta_0 + d\ff$, $\ff$ being a fiber of $\Sigma_d \to \PP^1$. Then the branch locus $\Delta_{f'}$ of 
$f'$ is of the form $\Delta_0 + \mcT$, where $\mcT$ is a reduced curve such that
(a) $\mcT$ is linear equivalent to $3\Delta$ and (b)
its singularities are at most simple. Under these
circumstances, $S$ is the canonical resolution of $S_O'$.
Choose a section $s$ such that $s\cdot O = 0$ and let 
$\sigma = t_{P_s}\circ[-1]$.

\begin{lem}\label{lem:fixed_locus_s}{The fixed locus $S^{\sigma}$ of $\sigma$ is of the form $R+ \Xi$, where
\begin{enumerate}
    \item[(i)] $R$ is horizontal and $R\cdot F = 4$,  
    \item[(ii)] $\Xi$ is fibral and
    \item[(iii)] $R\cdot s = R\cdot O = 0$.
\end{enumerate}
Here we say that a divisor $D$ is horizontal if no irreducible component of $D$ is contained in any fiber,
while $D$ is fibral if every irreducible component
is contained in a fiber.
}
\end{lem}

\proof By Lemma\ref{lem:fixed-points}, $S^{\sigma}$ is a reduced divisor on $S$. Put $S^{\sigma} = R + \Xi$, 
where $R$ is horizontaland $\Xi$ is fibral. Since $F/\langle\sigma\rangle \cong \PP^1$, $R\neq \emptyset$ and $R\cdot F = 4$. As $s\cdot O = 0$ and $s = \sigma(O)$, both $O$ and $s$ contain no fixed point. Hence $R\cdot O = R\cdot s = 0$.
 \endproof

 By Lemma\ref{lem:fixed_locus_s},
 the image 
 $f_{O, P_s}(s)$ is a section of $S/\langle\sigma\rangle \to \PP^1$ with $(f_{O, P_s}(O))^2 = -d/2$. By blowing down irreducible components in
 fibers of $S/\langle\sigma\rangle \to \PP^1$ not
 meeting $f_{O, P_s}(O)$ suitably, we see that the resulting
 relatively minimal model of $S/\langle\sigma\rangle$ is $\Sigma_{d/2}$, the Hirzebruch surface of degree $d/2$.
 In the diagram in the Introduction, $q_{O,P_s}$ is nothing but $S/\langle\sigma\rangle \to \Sigma_{d/2}$ and $\tilde{f}_{O,P_s}(S^{\sigma})$ is a reduced divisor $\tilde\mcQ$ on $\Sigma_{d/2}$ such that $\tilde\mcQ\cdot\ff = 4$. By the observation in Section \ref{sec:action}, $\tilde\mcQ$ has only simple singularities.  Hence $S$ is the canonical resolution of the double cover
 of $\Sigma_{d/2}$ with branch locus $\tilde\mcQ$.
 

 \subsection{Proof of Theorem~\ref{thm:concurrent}}
 
Let $\mcQ$ be an irreducible quartic with at most nodes. In this case, $\mcQ$ has  bitangent lines. Choose $z_o \in \PP^2\setminus \mcQ$.
As we explain in the Introduction, we have a rational elliptic surface
$\varphi_{\mcQ, z_o}: S_{\mcQ, z_o} \to \PP^1$ with
a split model $f_{O^+, O^-}: S'_{O^+, O^-} \to \Sigma_1$.  Since the elliptic fibration
$\varphi_{\mcQ, z_o} : S_{\mcQ, z_o} \to \PP^1$ is induced by the pencil of lines $\Lambda_{z_o}$ centered at $z_o$, singular fibers of $\varphi$ are 
given by members of $\Lambda_{z_o}$ not intersecting with $\mcQ$ transversely and types of singular fibers depend on how they intersect
with $\mcQ$. For such lines in $\Lambda_{z_o}$ not passing through
nodes of $\mcQ$, { i.e., tangent lines, bitangent lines and so on},
we refer to \cite[p.21, Table 1]{kuwata2005}. As for  lines
in $\Lambda_{z_o}$ through nodes of $\mcQ$, by our
observation in Section~\ref{sec:action}, we have the following table:

\begin{center}
\begin{tabular}{|c|c|} \hline
Type of fiber &  $\mcQ$ and a line $L \in \Lambda_{z_o}$  \\ \hline
$\I_2$ & $L$ passes a node and meet two other distinct points. \\
$\I_3$ & $L$ passes a node and tangent at other smooth points. \\
$\I_4$ & $L$ connects two nodes. \\
$\III$ & $L$ is tangent to one of smooth branches of the node. \\
$\IV$  & $L$ is inflectional tangent to one of smooth branch nodes. \\ \hline
\end{tabular}
\end{center}

Let $\mcQ$ be an irreducible quartic with $\alpha$
nodes as its singularities. Assume that there exist
$m$ bitangent lines $l_1, \ldots, l_m$ concurrent
at $z_o$. Then $\varphi_{\mcQ, z_o}: S_{\mcQ, z_o} \to \PP^1$ has at 
least $m$ singular fibers whose types are either $\I_2$ or $\III$.
We write their irreducible decomposition by
\[
F_{i} = \Theta_{i, 0} + \Theta_{i, 1}, \quad i = 1, \ldots, m.
\]
We may also assume that
\[
O^+\cdot \Theta_{i, 0} = 1 \quad O^-\cdot \Theta_{i,1} = 1, \quad i = 1, \ldots, m.
\]
Now we consider the ramified model with respect to
$O^+$,  $(E_S, O^+)$ 
We consider $O^-$ as an element in $E_{S}(\CC(t))$, $O^+$ being the zero, 
and we compute the height pairing $\langle O^-, O^-\rangle$. As $O^+\cdot O^- = 0$, we have

\[
0 \le \langle O^-, O^-\rangle  =  2 - \sum_{v \in \Red(\varphi)}{\mathrm{Contr}}_v(O^-) \le  2 - \frac m2.
\]
Hence $m \le 4$. 
Moreover, when $\mcQ$ has $3$ nodes, the configuration of singular fibers of 
$\varphi_{\mcQ, z_o}: S_{\mcQ, z_o} \to \PP^1$ contains at least $k$ $\I_2$ type, $l$ $\III$ type (k+ l = m) from $l_1, \ldots, l_m$ and
$3$ others from $3$ nodes. Since the sum of the topological Euler numbers is equal to $12$, we have $2k + 3l \le 6$.
Hence $m \le 3$. Now we consider possible pairs $(k, l)$ for the maximal case. We first note that $O^-$ is a $2$-torsion when
$m = 4$.

\underline{The case $(m, \alpha) = (4, 0)$.}  In this case, we have $2k + 3l \le 12$. If $(k, l) = (1, 3)$, 
we infer that the configuration of singular fibers of $\varphi$ is $\I_1, \I_2, 3\III$. On the other hand, there
exist no rational elliptic surface with a $2$-torsion having such a configuration of singular fibers by \cite{miranda90} or \cite{persson90}. Hence, $(1,3)$
does not occur. 

\underline{The case $(m, \alpha) = (4, 1)$.}  In this case , $2k + 3l \le 10$. Hence possible
$(k,l)$'s are $(4, 0), (3, 1)$ and $(2, 2)$.

\underline{The case $(m, \alpha) = (4, 2)$.}  In this case , $2k + 3l \le 8$. Hence only
$(k,l)= (4, 0)$.

\underline{The case $(m, \alpha) = (3, 3)$.} In this case, $2k + 3l \le 6$. Hence only
$(k, l) = (3, 0)$ is possible.

In next section, we show all the cases for $(k, l)$ occur by giving explicit examples.


%% file: example.tex
\section{Examples}\label{sec:examples}

In this section, we give examples for all the possible maximal cases in 
Theorem~\ref{thm:concurrent}.
Before we go on to give explict examples, we fix our settings after
 \cite{bannai-tokunaga15,bannai-tokunaga17}.
 Let $\Sigma_2$ be
  the Hirzebruch surface of degree $2$. Take affine open subsets $U_1$ and $U_2$ of $\Sigma_2$ so that
  \begin{enumerate}
  \item[(i)] $U_i = \CC^2$ $(i = 1, 2)$ and 
  \item[(ii)] the coordinates of $U_1$ and $U_2$ are $(t, x)$ and $(s, x')$, respectively with $s = 1/t$ and
  $x' = x/t^2$.
  \end{enumerate}

  Under these coordinates, the unique section $\Delta_0$ with $\Delta_0^2 =
  -2$ is given by $x = x' = \infty$ and $S'_O$ is obtained in the following way:
  
  
   There exists a reduced divisor $\mcT_S$ on $\Sigma_2$ such that (a) $\mcT_S$ is given by the equation of the form
  \[
  \mcT_S: F_{\mcT_S}(t, x) = x^3 + a_2(t)x^2 + a_{4}(t)x + a_{6}(t) = 0,
  \]
  where $a_{2i}(t) \in \CC[t], \deg a_{2i} \le di$ ($i = 1, 2, 3$) on  $U_1$ and 
  $O$ is given by the preimage of $\Delta_o$, (b) $\mcT_S$ has at most simple singularities
  (see \cite{bhpv} for simple singularities) and  (c) $(E_S, O)$ is given by the Weierstrass equation
  $y^2 = F_{\mcT_S}(t, x)$.


\subsection{The case \texorpdfstring{$\mathcal{Q}$}{Q} is smooth}

\begin{ex}\label{eg:ex1}

Consider a curve $\mcT^a$ on $U_1$ given by 
\[
F_1(t, x)= x(x^2 - 2(t^2+1)x -t^3 - 3t^2  - 2t) = 0.
\]
$\mcT^a$ can be extended  a reduced divisor $\mcT$ on $\Sigma_2$ satisfying the
conditions (a) and (b) on $\mathcal{T}_S$. Let $S_1$ be the canonical
resolutions of  the double cover, $S'_{O}$, whose generic fiber $(E_{S}, O)$ is given by
$y^2 = F_1(t, x)$.
Since $\mcT_1$ has nodes over $t=-2,-1,0,\infty$, the elliptic fibration
$\varphi: S \to \PP^1$ has $4$ singular fibers of type $\I_2$.
Let $P_0$ be a rational point $(0,0) \in E_S(\CC(t))$. We see that 
$\gamma(P_0) = [1,1,1,1]$, where $\gamma$ is the homomorphism given in Section~\ref{sec:notation}.
Here we list singular fibers over $t = -2, -1, 0, \infty$ in this order.
By applying (\ref{eq:transformation2}) to $P_0$, $w_{P_0}$ is explicitly given as follows:
\[
x = y'+x'^2 {+ t^2 + 1}, \quad
y = \sqrt{2}x'x.
\]
Then, we obtain a split model $E_{S, O, P_0}:y'^2=F_1'(t,x'):=(x'^2+t^2+1)^2 + t(t+1)(t+2)$.

The branch locus of $S'_{O,P_0} \to \PP^2$ is a smooth quartic $\mathcal{Q}$ given by $F_1'(t,x)=0$. 
In particular, $\mathcal{Q}$ is a reduced and smooth curve.
{By \cite[p. 21, Table 1]{kuwata2005}}, there are $4$ concurrent ordinary bitangent lines $t=-2,-1,0,\infty$.
Thus we have an example such that $\mathcal{Q}$  is smooth and 
$(k,l)=(4,0)$ in Theorem~\ref{thm:concurrent}.

\end{ex}


For the remaining cases $(k, l)  = (3, 1), (2, 2)$ and $(0,4)$, we consider the following equations and $P_0$:

\begin{center}
\begin{tabular}{|c|c|c|} \hline
$(k, l)$ & Equations of $\mcT$ & $P_0$ \\ \hline
$(3,1)$ & $F_2:= x(x^2 - 2x -t^3 - 3t^2  - 2t)$ & $(0,0)$ \\ \hline
$(2, 2)$ & $F_3:=x(x^2 - 2tx -t^3 - 3t^2  - 2t)$ & $(0,0)$ \\ \hline
$(0, 4)$ & $F_4:=x(x^2 -t^3 - 3t^2  - 2t)$ & $(0,0)$ \\ \hline
\end{tabular}
\end{center}

For each case, the corresponding elliptic surface $\varphi: S \to \PP^1$ has $4$ singular fibers over $t = -2, -1, 0$ and $\infty$.  Their types are follows:

\begin{center}
\begin{tabular}{|c|c|c|c|c|} \hline
        & $t=-2$ & $t = -1$ & $t=0$ & $t=\infty$ \\ \hline
$(3,1)$ & $\I_2$ & $\I_2$ & $\I_2$ & $\III$ \\ \hline
$(2, 2)$ & $\I_2$ & $\I_2$ & $\III$ & $\III$ \\ \hline
$(0, 4)$ & $\III$ & $\III$ & $\III$ & $\III$ \\ \hline
\end{tabular}
\end{center}

Now we repeat the same argument as the case of $(4, 0)$. We then obtain quartic curves $\mcQ$ with the desired bitangent lines 
as follows:

\begin{center}
\begin{tabular}{|c|c|c|c|} \hline
$(k, l)$ & Equations of $\mcQ$ & Ordinary bitangent  & Special bitangent  \\ \hline
$(3,1)$ & $F'_2:= (x'^2+1)^2 + t(t+1)(t+2)$ & $t=-2, -1, 0$ & $t = \infty$ \\ \hline
$(2, 2)$ & $F'_3:=(x'^2+ t)^2 + t(t+1)(t+2)$ & $t = -2, -1$ & $t = 0,  \infty$ \\ \hline
$(0, 4)$ & $F'_4:=(x'^2)^2 + t(t+1)(t+2)$ & None  & $t = -2, -1, 0, \infty$ \\ \hline
\end{tabular}
\end{center}
{Here we make use of \cite[p.21, Table 1]{kuwata2005}}.

\subsection{The case \texorpdfstring{$\mathcal{Q}$}{Q} has 1 node}

\begin{ex}\label{eg:ex5}
Consider a curve $\mcT^a$ on $U_1$ given by
\[
F_5(t,x):= (x-t^2)(x^2 -10tx +25x -36) =0
\]
$\mcT^a$ can be extended  a reduced divisor $\mcT$ on $\Sigma_2$ satisfying the
conditions (a) and (b) on $\mathcal{T}_S$. Let $S$ be the canonical
resolutions of  the double cover, $S'_{O}$, whose generic fiber $(E_{S}, O)$ is given by
$y^2 = F_5(t, x)$.
 $\mcT$ has nodes over $t=-1, 2, 3, 6,\infty$. Hence we see that
 the elliptic fibration
$\varphi: S \to \PP^1$ has $5$ singular fibers of type $\I_2$ over
 $t=-1, 2, 3, 6,\infty$.
Let $P_0 = (t^2,0) \in E_S(\CC(t))$. We see that 
$\gamma(P_0) = [1,1,1,1, 0]$, where $\gamma$ is the homomorphism given in Section~\ref{sec:notation}.
Here we list singular fibers over $t = -1, 2, 3, 6, \infty$ in this order.
By applying (\ref{eq:transformation2}) to $P_0$, the translation $w_{P_0}$ is explicitly given as follows:
\begin{eqnarray*}
x = y'+x'^2+\frac{5}{2}(2t-5) \qquad
y = \sqrt{2}x'(x-t^2).
\end{eqnarray*}
Then, we obtain a split model $E_{S, O, P_0}:y'^2=F_5'(t,x')$:
\[
F_5'(t,x') := \left(x'^2 - t^2 + 5t - \frac{25}{2} \right)^2 - (t + 1)(t - 2)(t - 3)(t - 6).
\]

The branch locus of $S'_{O,P_0} \to \PP^2$ is a quartic with one node $\mathcal{Q}$ given by $F_5'(t,x)=0$. 
{By \cite[p.21, Table 1]{kuwata2005}} and the second table in Section~\ref{sec:action}, 
there are $4$ concurrent ordinary bitangent lines $t=-1,2, 3, 6$. 
Thus we have an example for the case when
$\mathcal{Q}$ has one node and 
$(k,l)=(4,0)$ in Theorem~\ref{thm:concurrent}.

\end{ex}

\begin{ex}\label{eg:ex6}
Consider a curve $\mcT^a$ on $U_1$ given by
\[
F_6(t,x):= (x^2 -t^3 - 2t^2)(x-1) =0
\]
$\mcT^a$ can be extended  to a reduced divisor $\mcT$ on $\Sigma_2$ satisfying the
conditions (a) and (b) on $\mathcal{T}_S$. Let $S$ be the canonical
resolutions of  the elliptic surface $S'_{O}$ whose generic fiber $(E_{S}, O)$ is given by
$y^2 = F_5(t, x)$.
 $\mcT$ has $5$ nodes over $t=\frac{-1-\sqrt{5}}{2},-1,0,\frac{-1+\sqrt{5}}{2},\infty$. Hence we see that $\varphi : S \to \PP^1$ has $4$
 singular fibers of type $\I_2$ over $t=\frac{-1-\sqrt{5}}{2},-1,0,\frac{-1+\sqrt{5}}{2}$ and one of type $\III$ over $\infty$.
Choose a rational point  $P_0=(1,0) \in E_{S}(\CC(t))$. If we list singular fibers over
$t = \frac{-1-\sqrt{5}}{2},-1,0,\frac{-1+\sqrt{5}}{2},\infty$ in
this order, we have $\gamma(P_0) =[1,1,0,1,1]$.
Now, by applying (\ref{eq:transformation2}) to $P_0$,  $w_{P_0}$ is explicitly given as follows:
\begin{eqnarray*}
x = y'+x'^2 \qquad
y = \sqrt{2}x'(x-1).
\end{eqnarray*}
Then, we obtain a split model $S'_{O, P_0}:y'^2=F_6'(t,x')$:
\[
F_6'(t,x')= \left(x'^2-1 \right)^2 + (t+1)\left( t-\frac{-1-\sqrt{5}}{2} \right) \left(t-\frac{-1+\sqrt{5}}{2} \right).
\]
The branch locus of $S'_{O, P_0} \to \PP^2$ is another quartic $\mathcal{Q}$ given by $F_6'(t,x)=0$. 
In particular, $\mathcal{Q}$ is a reduced curve and has $1$ node at $(t,x')=(0,0)$. {By \cite[p.21, Table 1]{kuwata2005}} and the second table in Section~\ref{sec:action}, 
there are $4$ concurrent bitangent lines $t=\frac{-1-\sqrt{5}}{2},-1,\frac{-1+\sqrt{5}}{2},\infty$.
In particular, the bitangent line $t = \infty$ is special.
As a result, this case is that $\mathcal{Q}$ has $1$ node where $(k,l)=(3,1)$ in Theorem~\ref{thm:concurrent}.

\end{ex}

\begin{ex}\label{eg:ex7}
Consider a curve $\mcT^a$ in $U_1$ given by
\[
F_6(t,x):=(x^2 - t^3 - t^2) \left(x - \frac{3}{2}t - \frac{3}{2} \right) = 0.
\]
Again $\mcT^a$ can be extended to a reduced divisor $\mcT$ on $\Sigma_2$ satisfying the conditions  (a) and (b) on $\mcT_S$. Let $S$ be the canonical resolution of the 
elliptic surface $S'_O$ whose generic fiber $(E_S, O)$ is given by $y^2 = F_6(t, x)$. 
Note that this example is the one given in \cite{bty20}.
Since $\mcT$ has $5$ nodes over  $t=-1,-\frac{3}{4},0,3,\infty$, the elliptic fibration
has $5$ singular fibers of type $\I_2$ over $t=-\frac{3}{4},0,3$ and type $\III$
over $t = -1, \infty$.
Choose a rational point $P_0=(\frac{3}{2}t+\frac{3}{2},0) \in E_S(\CC(t))$. We denote
the corresponding section by $s_0$. Then we have
$\gamma(s_0) = [1,1,0,1,1]$. Here we list $5$ singular fibers over
$t=-1,-\frac{3}{4},0,3,\infty$ in this order.
Now, by applying (\ref{eq:transformation2}), the translation $w_{P_0}$ is explicitly given as follows:
\begin{eqnarray*}
x = y'+x'^2 \qquad
y = \sqrt{2}x' \left(x-\frac{3t+3}{2} \right).
\end{eqnarray*} 
Then, we obtain a split model $y'^2=F_6'(t,x'):= (x'^2 -\frac{3}{2}t - \frac{3}{2})^2+\frac{1}{4}(4t + 3)(t + 1)(t - 3)$.
{By \cite[p.21, Table 1]{kuwata2005}} and the second table in Section~\ref{sec:action},  the quartic $\mathcal{Q}$ given by $F_6'(t,x')=0$ is a reduced curve and has $1$ node at $(t,x')=(0,0)$.
There are $4$ concurrent bitangent lines $t=-1,-\frac{3}{4},3,\infty$ which are
ordinary for $t = -\frac{3}{4}, 3$ and special for $t = -1, \infty$.
Hence we have an example such that $\mathcal{Q}$ has $1$ node where $(k,l)=(2,2)$ in Theorem~\ref{thm:concurrent}.

\end{ex}

\subsection{The cases \texorpdfstring{$\mathcal{Q}$}{Q} have 2 and 3 nodes}

\begin{ex}
\label{ex:LLQ}

Consider a curve $\mcT^a$ in $U_1$ given by
\[
F_7(t, x):= (x-t^2)(x-3t+2)(x+3t+2) = 0.
\]
$\mcT^a$ can also be extended to a reduced divisor $\mcT$ on $\Sigma_2$ (a union of
three sections) satisfying the conditions (a) and (b) on $\mcT_S$.
Let $S$ be the canonical resolutions  of $S'_O$ whose generic fiber is given 
by $y^2 = F_7(t, x)$. Note that this example is the one  given in \cite[EXAMPLE 5.2]{tokunaga14}. Since $\mcT$ has $6$ nodes, we see that $\varphi: S \to \PP^1$ has $6$ singular fibers
over $t=-2,-1,0,1,2,\infty$ and all of them are of type $\I_2$. Choose rational
points $P_0=(3t-2,0)$ and $P_1=(t+2,2\sqrt{2}(t-2)(t+1))$ in $E_S(\CC(t))$.
If we list singular
fibers over $t=-2,-1,0,1,2,\infty$ in this order, 
we have $\gamma(P_0) = [0, 0, 1, 1, 1,1]$ and $\gamma(P_1) = [0, 1, 0, 0, 1,1]$.
%
%
By applying (\ref{eq:transformation2}) to $P_0$,  $w_{P_0}$ is explicitly given as follows:
\begin{eqnarray*}
x = y'+x'^2+\frac{1}{2}(t^2-3t-2) \qquad
y = \sqrt{2}x'(x-3t+2).
\end{eqnarray*}
Then, we obtain a split model $E_{S, O, P_0} :y'^2=F_7'(t,x'): =(x'^2 + \frac{1}{2}t^2 - \frac{9}{2}t + 1)^2 + 6(t - 1)(t - 2)t$.
The branch locus of $S'_{O, P_0}$ is a quartic $\mathcal{Q}$ given by $F_7'(t,x')=0$.
In particular, $\mathcal{Q}$ is a reduced and has $2$ node at $(t,x')=(-2, 0)$, $(-1,0)$.
{By \cite[p.21, Table 1]{kuwata2005}} and the second table in Section~\ref{sec:action},  there are $4$ ordinary concurrent bitangent lines $t=0,1,2,\infty$.
Hence we have an example for the case $2$ nodes where $(k,l)=(4,0)$ in Theorem~\ref{thm:concurrent}.

We next apply (\ref{eq:transformation2}) to $P_1$. Then $w_{P_1}$ is explicitly given
as follows:
\begin{eqnarray*}
x = y'+x'^2+\frac{1}{2}(t^2-t-6) \qquad
y = \sqrt{2}x'(x-t-2) - 2\sqrt{2}(t-2)(t+1).
\end{eqnarray*}
Then, we obtain a split model $y'^2=F_8'(t,x'):=(x'^2 + \frac{1}{2}t^2 - \frac{3}{2}t - 5)^2 + 2(t - 4x' + 8)(t + 1)(t - 2)$.
The quartic $\mathcal{Q}$ given by $F_8'(t,x')=0$ is a reduced curve and has $3$ node at $(t,x')=(-2,-2)$, $(0,-1)$, $(1,-2)$.
{By \cite[p.21, Table 1]{kuwata2005}} and the second table in Section~\ref{sec:action},  there are $3$ ordinary concurrent bitangent lines $t=-1,2,\infty$.
This example gives the case for $3$ nodes where $(k,l)=(3,0)$ in Theorem~\ref{thm:concurrent}.
\end{ex}


